\documentclass{article}
\usepackage{latexsym,amsmath,amssymb,graphicx}
\input amssym.def \input amssym
\newtheorem{te}{Theorem}[section]
\newtheorem{de}[te]{Definition}
\newtheorem{lm}[te]{Lemma}
\newtheorem{fa}[te]{Fact}
\newtheorem{pp}[te]{Proposition}
\newtheorem{co}[te]{Corollary}
\newtheorem{ex}[te]{Example}
\def\dokaz{\noindent{\bf Proof. }}
\def\kraj{\hfill $\Box$ \par \vspace*{2mm} }

\def\widemid{\hspace{1mm}\widetilde{\mid}\hspace{1mm}}

\def\po{\exists}
\def\str{\rightarrow}
\def\Str{\Rightarrow}
\def\rtS{\Leftarrow}
\def\dl{\Leftrightarrow}

\def\ps{\subseteq}

\def\ra{\setminus}

\def\rest{\upharpoonright}

\def\cU{{\cal U}}
\def\cA{{\cal A}}

\begin{document}
\begin{center}
           {\large \bf Divisibility in the Stone-\v Cech compactification}
\end{center}
\begin{center}
{\small \bf Boris  \v Sobot}\\[2mm]
{\small  Department of Mathematics and Informatics, University of Novi Sad,\\
Trg Dositeja Obradovi\'ca 4, 21000 Novi Sad, Serbia\\
e-mail: sobot@dmi.uns.ac.rs}
\end{center}
\begin{abstract} \noindent
After defining continuous extensions of binary relations on the set $N$ of natural numbers to its Stone-\v Cech compactification $\beta N$, we establish some results about one of such extensions. This provides us with one possible divisibility relation on $\beta N$, $\widemid$, and we introduce a few more, defined in a natural way. For some of them we find equivalent conditions for divisibility. Finally, we mention a few facts about prime and irreducible elements of $(\beta N,\cdot)$. The motivation behind all this is to try to translate problems in number theory into $\beta N$.
\vspace{1mm}\\

{\sl 2010 Mathematics Subject Classification}:
54D35, %Extensions of spaces
54D80, %Special constructions of spaces
22A15. %Structure of topological semigroups

{\sl Key words and phrases}: divisibility, Stone-\v Cech compactification, ultrafilter
\end{abstract}

\section{Introduction}

It is well-known that every semigroup $(S,\cdot)$, supplied with the discrete topology, can be extended to a compact Hausdorff right-topological semigroup $(\beta S,\odot)$. $\beta S$, {\it the Stone-\v Cech compactification} of the space $S$, is the space of all ultrafilters over $S$. The topology on $\beta S$ is generated by (clopen) base sets of the form ${\bar A}=\{p\in\beta S:A\in p\}$. We will do this with the semigroup $(N,\cdot)$, where $\cdot$ is the usual multiplication.

Let us first fix some notation: for $A\ps N$ we denote $A^c=N\setminus A$, and for infinite $A\ps N$, $[A]^{\aleph_0}$ is the set of all infinite subsets of $A$.

For $A\ps N$ and $n\in N$ let $A/n=\{m\in N:mn\in A\}=\{\frac an:a\in A,n\mid a\}$. We will frequently use the following fact.

\begin{fa}\label{presek}
$A/n\cap B/n=(A\cap B)/n$.
\end{fa}

The extended operation on $\beta N$ will also be denoted by $\cdot$ and is defined in the following way:
$$A\in p\cdot q\dl\{n\in N:A/n\in q\}\in p.$$
This is the unique way of defining $\cdot$ so that the left-multiplication functions $\lambda_n:x\mapsto nx$ are continuous for $n\in N$, and all right-multiplication functions $\rho_p:x\mapsto xp$ are continuous for $p\in\beta N$. For each $n\in N$ the principal ultrafilter $p_n=\{A\ps N:n\in A\}$ is identified with the respective element $n$. The singletons $\overline{\{n\}}=\{n\}$ are clopen. If $A\in [N]^{\aleph_0}$ it is also common to think of $\beta A$ as a subspace of $\beta N$, and to denote $A^*=\beta A\ra A$. Clopen sets in $\beta N$ are exactly those of the form $\bar A$, and clopen sets in $N^*$ are exactly $A^*$ for $A\subseteq N$.

The book \cite{HS} is a rich supply of facts about algebra in the Stone-\v Cech compactification and we will use results from it heavily in this paper. We begin with some basic facts. The topological center of a right-topological semigroup $(S,\cdot)$ is the set $\{p\in S:\lambda_p\mbox{ is continuous}\}$, and the algebraical center is $\{p\in S:\forall x\in S\;px=xp\}$.

\begin{pp}\label{centar}
(\cite{HS}, Theorem 6.10) $N$ is both the algebraic and the topological center of $(\beta N,\cdot)$.
\end{pp}

\begin{pp}\label{Gdelta}
(\cite{HS}, Theorem 3.36) Every nonempty $G_\delta$ subset of $N^*$ has a nonempty interior in $N^*$.
\end{pp}

(This is just a restatement of the fact that the pseudointersection number ${\goth p}$ is greater than $\aleph_0$.)

\section{Extending functions and relations}\label{extending}

Since we supplied $N$ with the discrete topology, every function from $N$ to any topological space is continuous. Perhaps the most important property of the Stone-\v Cech compactification is the following.

\begin{pp}\label{extfun}
(\cite{HS}, Theorem 3.27) If $Y$ is any compact Hausdorff topological space, then every function $f:N\str Y$ extends in a unique way to a continuous function $\widetilde{f}:\beta N\str Y$.
\end{pp}

In particular, every function $f:N\str N$ extends uniquely to continuous $\widetilde{f}:\beta N\str\beta N$. In this case we have the following.

\begin{pp}\label{extoblik}
(\cite{HS}, Lemma 3.30) If $f:N\str N$, then: (a) for every $p\in\beta N$, $\widetilde{f}(p)=\{A\ps N:f^{-1}[A]\in p\}$; (b) for all $A\in p$, $f[A]\in\widetilde{f}(p)$.
\end{pp}

\begin{pp}\label{exthom}
(\cite{HS}, Corollary 4.22) Let $(Y,\cdot)$ be a compact Hausdorff right-topological semigroup. If $f:N\str Y$ is a homomorphism such that $f[N]$ is contained in the topological center of $(Y,\cdot)$, then $\widetilde{f}:\beta N\str Y$ is also a homomorphism.
\end{pp}

\begin{pp}\label{extbij}
(\cite{HS}, Exercise 3.4.1) Let $D$ and $E$ be discrete spaces. (a) If $f:D\str E$ is injective, so is $\widetilde{f}:\beta D\str\beta E$; (b) if $f:D\str E$ is surjective, so is $\widetilde{f}:\beta D\str\beta E$.
\end{pp}

If $f:N\str N$ is injective, then $\widetilde{f}$ is a bijection between $\beta N$ and $\overline{f[N]}$. For example, let $f(x)=x^2$. Then, since $f[N]=\{x^2:x\in N\}$, $\widetilde{f}$ maps every $p\in\beta N$ to the ultrafilter $\widetilde{f}(p)$ generated by $\{\{x^2:x\in A\}:A\in p\}$.\\

We would like to be able to extend arbitrary binary relations in this way (we will work only with $\beta N$, but the results here are actually quite general). To make this more precise, we fix some natural notation. We say that a relation $\sigma\ps(\beta N)^2$ is an {\it extension} of $\rho\ps N^2$ if $\sigma\rest_{N^2}=\rho$. If $A\ps X$ and $\rho\subseteq X^2$, let $\rho[A]=\{n\in N:\exists m\in A\;(m,n)\in\rho\}$ be the direct image, and $\rho^{-1}[A]=\{m\in N:\exists n\in A\;(m,n)\in\rho\}$ the inverse image of $A$ by $\rho$.

\begin{de}
If $X$ is a topological space, a binary relation $\sigma\ps X^2$ is {\it continuous} if
\begin{equation}\label{ROP}\tag{ROP}
\mbox{for every open set }U\ps X\mbox{ the set }\sigma^{-1}[U]\mbox{ is also open.}
\end{equation}
\end{de}

\begin{ex}
Unlike functions, for binary relations the condition
\begin{equation}\label{RCL}\tag{RCL}
\mbox{for every }q\in X\mbox{ the set }\sigma^{-1}[\{q\}]\mbox{ is closed}
\end{equation}
does not follow from (\ref{ROP}). To see this, let $\sigma=N\times(\beta N)$. It is continuous because every subset of $N$ is open. But for any $q\in\beta N$ the set $\sigma^{-1}[\{q\}]=N$ is not closed.
\end{ex}

By analogy with the extension $\widetilde{f}$ from Proposition \ref{extoblik} we define an extension of a binary relation.

\begin{de}\label{deftilda}
For $\rho\ps N^2$ let $\widetilde{\rho}\ps(\beta N)^2$ be given by $\widetilde{\rho}=\{(p,q)\in(\beta N)^2:\forall A\in p\;\rho[A]\in q\}$.
\end{de}

\begin{te}\label{rhotilda}
Let $\rho\ps N^2$ be a binary relation. Then $\widetilde{\rho}$ is a continuous extension of $\rho$ satisfying (\ref{RCL}).
\end{te}

\dokaz For $m,n\in N$ we have $m\widetilde{\rho}n$ if and only if $n\in\rho[\{m\}]$, i.e.\ if and only if $m\rho n$, so $\widetilde{\rho}$ is an extension of $\rho$.

Now we prove (\ref{ROP}). Let $\bar B$ be a base set and $p\in\widetilde{\rho}^{-1}[\bar B]$. So let $q\in\beta N$ be such that $p\widetilde{\rho}q$, and $B\in q$. We define a function $f:N\str N$ as follows: $f(x)=\min\{y\in B:x\rho y\}$ if such $y$ exists, and $f(x)\in B$ is arbitrary otherwise. Obviously, $\widetilde{f}(p)\in{\bar B}$. Since $\widetilde{f}$ is continuous, there is $A\in p$ such that for all $r\in{\bar A}$ we have $\widetilde{f}(r)\in{\bar B}$. The set $\rho^{-1}[B]$ is in $p$ (otherwise its complement $(\rho^{-1}[B])^c$ would be in $p$, so by definition of $\widetilde{\rho}$ we would have $\rho[(\rho^{-1}[B])^c]\in q$, which is impossible since $\rho[(\rho^{-1}[B])^c]\cap B=\emptyset$). So, without loss of generality, we can assume that $A\ps\rho^{-1}[B]$. This means that for every $r\in{\bar A}$ holds $f[X\cap A]\ps\rho[X\cap A]$ for all $X\in r$. But $\widetilde{f}(r)$ is the unique element in $\bigcap\{f[X]:X\in r\}=\bigcap\{f[X\cap A]:X\in r\}$, and similarly $\widetilde{\rho}[\{r\}]=\bigcap\{\rho[X\cap A]:X\in r\}$, so $\widetilde{f}(r)\in\widetilde{\rho}[\{r\}]$ which means that $r\in\widetilde{\rho}^{-1}[{\bar B}]$. This means that the set $\widetilde{\rho}^{-1}[{\bar B}]$ is open.

To prove (\ref{RCL}), let $q\in\beta N$ and $p\notin\widetilde{\rho}^{-1}[\{q\}]$. Then there is $A\in p$ such that $\rho[A]\notin q$. It follows that $\bar A$ is the neighborhood of $p$ disjoint with $\widetilde{\rho}^{-1}[\{q\}]$.\kraj

\begin{ex}\label{notunique}
This extension is not unique. For example, if $n\in N$ and $q\in N^*\ra\widetilde{\rho}[\{n\}]$, let us show that the relation $\sigma=\widetilde{\rho}\cup\{(n,q)\}$ is also continuous. (\ref{ROP}) is obvious for open sets $U$ not containing $q$. On the other hand, for any $B\in q$ we have an open neighborhood $\{n\}$ of $n$ such that $\{n\}\ps\sigma^{-1}[{\bar B}]$, and for $p\neq n$ the existence of such neighborhood is guaranteed by the condition (\ref{ROP}) for $\widetilde{\rho}$. Regarding (\ref{RCL}), we need only to show that $\sigma^{-1}[\{q\}]$ is closed. But $\sigma^{-1}[\{q\}]=\widetilde{\rho}^{-1}[\{q\}]\cup\{n\}$, which is a union of two closed sets.
\end{ex}

\begin{te}\label{jedinstvenost}
Let $\rho\ps N^2$ be a binary relation and let $\sigma\ps(\beta N)^2$ be a continuous relation such that $\sigma\rest_{N\times\beta N}=\widetilde{\rho}\rest_{N\times\beta N}$. Then $\sigma\ps\widetilde{\rho}$.
\end{te}

\dokaz Assume the opposite, that $(p,q)\in\sigma\ra\widetilde{\rho}$. $(p,q)\notin\widetilde{\rho}$ means that there is $B\in p$ such that $\rho[B]\notin q$, i.e.\ $q\in\overline{(\rho[B])^c}$. By (\ref{ROP}) there is a base set $\bar A$ such that $p\in\bar A\ps\sigma^{-1}\left[\overline{(\rho[B])^c}\right]$. We may also assume that $A\ps B$. Choose any $a\in A$. There is $r\in\overline{(\rho[B])^c}$ such that $(a,r)\in\sigma$. On the other hand we have $r\notin\overline{\rho[B]}$, so by the definition of $\widetilde{\rho}$, $(a,r)\notin\widetilde{\rho}$, a contradiction since $(a,r)\in N\times\beta N$.\kraj

The example \ref{notunique} shows that we needed to assume that $\sigma$ and $\widetilde{\rho}$ agree on $N\times\beta N$. Now we prove several more things about $\widetilde{\rho}$.

\begin{te}\label{rsttilda}
(a) $\widetilde{\rho^{-1}}=(\widetilde{\rho})^{-1}$.

(b) If $\rho$ is reflexive, then so is $\widetilde{\rho}$. If $\rho$ is symmetric, then so is $\widetilde{\rho}$. If $\rho$ is transitive, then so is $\widetilde{\rho}$.
\end{te}

\dokaz (a) Assume that $(p,q)\in(\widetilde{\rho})^{-1}$ but $(p,q)\notin\widetilde{\rho^{-1}}$. Then there is $A\in p$ such that $\rho^{-1}[A]\notin q$. This means that $(\rho^{-1}[A])^c\in q$ so, since $(q,p)\in\widetilde{\rho}$, $\rho[(\rho^{-1}[A])^c]\in p$. But $A\cap\rho[(\rho^{-1}[A])^c]=\emptyset$, a contradiction.

We proved that $(\widetilde{\rho})^{-1}\ps\widetilde{\rho^{-1}}$, so by exchanging the roles of $\rho$ and $\rho^{-1}$ we have $(\widetilde{\rho^{-1}})^{-1}\ps\widetilde{\rho}$, which is equivalent to $\widetilde{\rho^{-1}}\ps(\widetilde{\rho})^{-1}$.

(b) For symmetry the statement follows from (a).

To see that $\widetilde{\rho}$ is reflexive, note that by reflexivity of $\rho$ we have $\rho[A]\supseteq A$ for all $A\ps N$. Hence $p\in\bigcap\{\overline{\rho[A]}:A\in p\}$, so $p\widetilde{\rho}p$.

Finally, let $p\widetilde{\rho}q$ and $q\widetilde{\rho}r$. For $A\in p$ we have $\rho[A]\in q$, which further implies $\rho[\rho[A]]\in r$. By transitivity of $\rho$, $\rho[\rho[A]]\ps\rho[A]$. Hence for all $A\in p$ we have $\rho[A]\in r$, so $p\widetilde{\rho}r$.\kraj

\begin{ex}
If $\rho$ is antisymmetric $\widetilde{\rho}$ need not be, and the same holds for irreflexivity. Indeed, it is not hard to see that, if $\leq_N$ and $<_N$ are the usual linear orders on $N$, then $\widetilde{\leq_N}=\leq_N\cup(\beta N\times N^*)$ and $\widetilde{<_N}=<_N\cup(\beta N\times N^*)$.
\end{ex}

If $h:X\str Y$ is a function, by $\ker h$ we denote its kernel relation on $X$: $(x_1,x_2)\in\ker h\dl h(x_1)=h(x_2)$.

\begin{te}\label{kernel}
If $h:N\str N$ and $\rho=\ker h$, then $\widetilde{\rho}=\ker\widetilde{h}$.
\end{te}

\dokaz Note that $(p,q)\in\ker\widetilde{h}$ if and only if for all $A\ps N$ holds $h^{-1}[A]\in p\dl h^{-1}[A]\in q$. In other words, $(p,q)\in\ker\widetilde{h}$ if and only if $p$ and $q$ agree on the unions of equivalence classes of $\ker h$.

We will prove that the equivalence class $\{q:(p,q)\in\ker\widetilde{h}\}$ of any $p\in\beta N$ is the same as $\widetilde{\rho}[\{p\}]=\{\overline{\rho[A]}:A\in p\}$. First, let $(p,q)\in\ker\widetilde{h}$ and $A\in p$. Then $\rho[A]$ is the union of all equivalence classes of $\rho=\ker h$ that have nonempty intersection with $A$. Since $A\ps\rho[A]$, it follows that $\rho[A]\in p$. Since $(p,q)\in\ker\widetilde{h}$, $\rho[A]\in q$.

Now let $q\in\widetilde{\rho}[\{p\}]$. For any set $A$ that is a union of equivalence classes of $\ker h$ we have $\rho[A]=A$, so for each such set $A\in p$ implies $A\in q$, and $A^c\in p$ implies $A^c\in q$. This means that $(p,q)\in\ker\widetilde{h}$.\kraj

\section{Ideals and cancelation}

A set $L\ps\beta N$ is {\it a left ideal} if $\beta N\cdot L\ps L$, i.e.\ if for all $p\in\beta N$ and all $x\in L$ we have $px\in L$. A left ideal $L$ is {\it minimal} if there is no $L_1\subset L$ that is also a left ideal. Right ideals and minimal right ideals are defined analogously. $L$ is {\it an ideal} if it is both left and right ideal.\\

Since the semigroup $(\beta N,\cdot)$ has a neutral element 1, each set of the form $\beta Np=\{xp:x\in\beta N\}$, where $p\in\beta N$, is a left ideal containing $p$, each set of the form $p\beta N=\{px:x\in\beta N\}$ is a right ideal containing $p$, and sets of the form $\beta Np\beta N=\{xpy:x,y\in\beta N\}$ are (two-sided) ideals. These are called principal left ideals, principal right ideals and principal ideals respectively. Clearly, each minimal left ideal $L$ is principal: $L=\beta Np$ for every $p\in L$.

\begin{pp}\label{ideal*}
(\cite{HS}, Theorem 4.36) $N^*$ is an ideal of $\beta N$.
\end{pp}

\begin{pp}\label{minideali}
(a) (\cite{HS}, Corollary 2.6, Theorem 2.7(a)) Every left ideal in $\beta N$ contains a minimal left ideal. Every right ideal in $\beta N$ contains a minimal right ideal.

(b) (\cite{HS}, Theorem 2.7(d)) If $L$ is a left ideal and $R$ is a right ideal then $L\cap R\neq\emptyset$.

(c) (\cite{HS}, Theorem 6.30, Corollary 6.41) There are $2^{\goth c}$ (disjoint) minimal left ideals of $\beta N$. There are $2^{\goth c}$ minimal right ideals of $\beta N$.
\end{pp}

\begin{pp}\label{idealK}
(a) (\cite{HS}, Theorem 1.51) $\beta N$ has the smallest ideal, denoted by $K(\beta N)$: an ideal contained in all other ideals of $\beta N$.

(b) (\cite{HS}, Theorem 1.64) $K(\beta N)=\bigcup\{L:L\mbox{ is a minimal left ideal of }\beta N\}=\bigcup\{R:R\mbox{ is a minimal right ideal of }\beta N\}$.

(c) (\cite{HS}, Theorem 2.10) The following conditions are equivalent: (i) $p\in K(\beta N)$ (ii) for all $q\in\beta N$ holds $p\in\beta Nqp$ (iii) for all $q\in\beta N$ holds $p\in pq\beta N$.
\end{pp}

By Proposition \ref{idealK}(c) we can not hope for a unique factorization result in $\beta N$. Hence cancelation laws are important for us in trying to establish how different factorizations of an element may be.

\begin{pp}\label{cancN}
(\cite{HS}, Lemma 8.1) If $n\in N$ and $p,q\in\beta N$, then $np=nq$ implies $p=q$.
\end{pp}

\begin{pp}\label{cancN2}
(\cite{HS}, Lemma 6.28) If $m,n\in N$ and $p\in\beta N$, then $mp=np$ implies $m=n$.
\end{pp}

\begin{pp}
(\cite{HS}, Theorem 8.7) $p\in\beta N$ is right cancelable if and only if for every $A\ps N$ there is $B\ps A$ such that $A=\{x\in N:B/x\in p\}$.
\end{pp}

\begin{pp}
(\cite{HS}, Theorem 8.10) The set of right cancelable elements contains a dense open subset of $N^*$, i.e.\ for every $U\in[N]^{\aleph_0}$ there is $V\in[U]^{\aleph_0}$ such that all $p\in{\bar V}$ are right cancelable.
\end{pp}

\section{Divisibility on $\beta N$}

There are several ways to extend the divisibility relation $\mid$ on $N$ to $\beta N$ that seem natural.

\begin{de}
Let $p,q\in\beta N$.

(a) $q$ is left-divisible by $p$, $p\mid_L q$, if there is $r\in\beta N$ such that $q=rp$.

(b) $q$ is right-divisible by $p$, $p\mid_R q$, if there is $r\in\beta N$ such that $q=pr$.

(c) $q$ is mid-divisible by $p$, $p\mid_M q$, if there are $r,s\in\beta N$ such that $q=rps$.

(d) $\widemid$ is the extension of $\mid$ as defined in Definition \ref{deftilda}.
\end{de}

It is obvious that the first three of these relations are reflexive and transitive. We can use Proposition \ref{idealK}(c) to show that neither of them is antisymmetric. For example, if $p\in K(\beta N)$, then for any $q\in\beta N$ we have $pq\mid_R p$ and, trivially, $p\mid_R pq$. From Proposition \ref{cancN2} follows that for $q\in N$ holds $pq\neq p$.\\

By \ref{rsttilda} $\widemid$ is also reflexive and transitive. From Theorem \ref{desnacont} it will follow that it is not antisymmetric either.\\

It is often useful to translate divisibility relations into statements about ideals.

\begin{fa}
(a) $p\mid_L q$ if and only if $\beta Nq\ps\beta Np$.

(b) $p\mid_R q$ if and only if $q\beta N\ps p\beta N$.

(c) $p\mid_M q$ if and only if $\beta Nq\beta N\ps\beta Np\beta N$.
\end{fa}

Clearly, $\mid_L\ps\mid_M$ and $\mid_R\ps\mid_M$. We will show that the following holds (all inclusions are strict):

$$\begin{array}{c}
\mid_L\\
\\
\mid_R
\end{array}\hspace{-3mm}
\begin{array}{c}
\rotatebox[origin=c]{-45}{$\subset$}\\
\rotatebox[origin=c]{45}{$\subset$}\\
\end{array}\hspace{-2mm}
\mid_M\subset\widemid$$

In general, much more is known about left ideals (and hence about left divisibility) than about right ones. However, for the right divisibility the following holds.

\begin{te}\label{desnacont}
The relation $\mid_R$ is a continuous extension of $\mid$ to $\beta N$.
\end{te}

\dokaz We need to prove that for each element $p\in\beta N$ and a neighborhood $\bar U$ of an element $q$ such that $p\mid_R q$ there is a neighborhood $\bar V$ of $p$ such that $\bar V\ps\mid_R^{-1}[{\bar U}]$. $p\mid_R q$ means that there is $r\in\beta N$ such that $q=pr=\rho_r(p)$. Since $\rho_r$ is continuous, there is $V\in p$ such that ${\bar V}\ps\rho_r^{-1}[{\bar U}]\ps\mid_R^{-1}[{\bar U}]$.\kraj

Now $\mid_R\subseteq\widemid$ follows from Theorems \ref{jedinstvenost} and \ref{desnacont}. $\mid_M\ps\widemid$ will be proven in Lemma \ref{odnosi}. In Lemma \ref{nijeRC2} we will prove that $\mid_M$ does not satisfy (\ref{RCL}), which will give us the strict inclusion $\mid_M\subset\widemid$ by Theorem \ref{rhotilda}.

\begin{lm}
No two of the relations $\mid_L$, $\mid_R$ and $\mid_M$ are the same. Neither of $\mid_L$ and $\mid_R$ is a subset of the other.
\end{lm}

\dokaz Let $p\in K(\beta N)$. Then, by Proposition \ref{idealK}(c), $q\mid_M p$ for all $q\in\beta N$. On the other hand, by Proposition \ref{minideali}(c) and Proposition \ref{idealK}(b) $K(\beta N)$ consists of $2^{\goth c}$ disjoint minimal left ideals; let $L_1$ and $L_2$ be two of them. If $l_1\in L_1$ and $l_2\in L_2$, we have $L_1=\beta Nl_1$ and $L_2=\beta Nl_2$. This means that no element of $\beta N$ can be left-divisible by both $l_1$ and $l_2$, say $l_1\nmid p$. So $\mid_L\subset\mid_M$.

Now let $R$ be a minimal right ideal. If $r$ is any of its elements, we have $R=r\beta N$. By Proposition \ref{minideali}(b) there must exist elements $r_1\in R\cap L_1$ and $r_2\in R\cap L_2$. Since $r_2\in R=r_1\beta N$, we have $r_1\mid_R r_2$, but since $r_2\notin L_1$, $r_1\mid_L r_2$ does not hold.

Analogously we conclude that $\mid_R\subset\mid_M$ and that $\mid_L\subseteq\mid_R$ does not hold.\kraj

\section{Divisibility by elements of $N$}

\begin{lm}
If $n\in N$, each of the statements: (i) $n\mid_L p$, (ii) $n\mid_R p$, (iii) $n\mid_M p$, (iv) $n\widemid p$ and (v) $nN\in p$ are equivalent.
\end{lm}

\dokaz Since elements $n\in N$ commute with all elements of $\beta N$ (by Proposition \ref{centar}), the first three statements are equivalent.

(ii)$\Str$(v) If $p=nq$, this means that $A\in p$ if and only if $\{s\in N:A/s\in q\}\in n$, i.e.\ if and only if $A/n\in q$. But $(nN)/n=N\in q$, so $nN\in p$.

(v)$\rtS$(ii) Let $nN\in p$ and $\cA=\{A/n:A\in p\}$. All elements of $\cA$ are nonempty, because $(nN)^c\notin p$. By Fact \ref{presek} $\cA$ has the finite intersection property, so there is an ultrafilter $q$ containing $\cA$, and $p=nq$.

Finally, $n\widemid p$ if and only if the direct image $\mid[\{n\}]=nN$ is in $p$.\kraj

In this case, when $n\in N$, we can drop the subscript and write only $n\mid p$.\\

It is clear from Proposition \ref{ideal*} that elements of $N$ can not be divisible by elements of $N^*$. On the other hand, an ultrafilter can be divisible by many natural numbers.

\begin{te}\label{zadatidelioci}
Let $A\ps N$ be downward closed for $\mid$ and closed for the operation of least common multiple. Then there is a nonempty $G_\delta$ set $X\ps N^*$ such that all $x\in X$ are divisible by all $n\in A$, and not divisible by any $n\notin A$.
\end{te}

\dokaz We want to show that the set $X=\bigcap_{n\in A}\overline{nN}\cap\bigcap_{n\notin A}\overline{(nN)^c}$ is nonempty. Since all the sets in the family $\cA=\{\overline{nN}:n\in A\}\cup\{\overline{(nN)^c}:n\notin A\}$ are closed, by compactness it suffices to show that $\cA$ has the finite intersection property. So let $a_1,a_2,\dots,a_k\in A$ and $b_1,b_2,\dots,b_l\notin A$ be given, and let $d$ be the least common multiple of $a_1,a_2,\dots,a_k$. Then $b_i\nmid d$ for $i=1,2,\dots,l$, because otherwise we would have $b_i\in A$. Hence $d\in\overline{a_iN}$ for $i=1,2,\dots,k$ and $d\in\overline{(b_iN)^c}$ for $i=1,2,\dots,l$.\kraj

\section{Equivalent conditions of divisibility}

If $p\in\beta N$, we denote $C(p)=\{A\ps N:\forall n\in N\;A/n\in p\}$ and $D(p)=\{A\ps N:\{n\in N:A/n=N\}\in p\}$. Using Fact \ref{presek} it is easy to see that both $C(p)$ and $D(p)$ are filters contained in $p$.

\begin{te}\label{Codp}
The following conditions are equivalent:

(i) $p\mid_L q$;

(ii) $C(p)\ps q$;

(iii) $C(p)\ps C(q)$.
\end{te}

\dokaz The equivalence of (i) and (ii) is Theorem 6.18 from \cite{HS}. (iii) implies (ii) is obvious because $C(q)\ps q$, so let us assume that $p\mid_L q$. If there would exist $A\in C(p)\setminus C(q)$, there would also exist an ultrafilter $r\supseteq C(q)\cup\{A^c\}$ so we would have $p\mid_L q$ and $q\mid_L r$, but not $p\mid_L r$, a contradiction.\kraj

Let $\cU=\{S\ps N:S\mbox{ is upward closed for }\mid\}$.

\begin{te}\label{Dodp}
The following conditions are equivalent:

(i) $p\widemid q$, i.e.\ for all $A\ps N$, $A\in p$ implies $\mid[A](=\{m\in N:\po a\in A\;a\mid m\})\in q$;

(ii) $p\cap\cU\ps q\cap\cU$;

(iii) $D(p)\ps D(q)$;

(iv) $D(p)\ps q$.
\end{te}

\dokaz (i)$\Str$(ii) Let $A\in p\cap\cU$. Then $A=\mid[A]$, so by (i) $A\in q$ as well.

(ii)$\Str$(iii) For any $A\ps N$ the set $B_A=\{n\in N:A/n=N\}=\{n\in N:nN\ps A\}$ is in $\cU$. So by (ii) $B_A\in p$ implies $B_A\in q$, which means that $A\in D(p)$ implies $A\in D(q)$.

(iii)$\Str$(iv) is obvious.

(iv)$\Str$(i) Let $A\in p$. For any $n\in A$ we have $nN\ps\mid[A]$. Hence the set $\{n\in N:nN\ps\mid[A]\}$ contains $A$, so it is in $p$. It follows that $\mid[A]\in D(p)$, thus $\mid[A]\in q$.\kraj

\begin{co}
Let $f:N\str N$ be a function.
(a) If $f(m)\mid m$ for all $m\in N$, then $\widetilde{f}(p)\widemid p$ for all $p\in\beta N$.

(b) If $m\mid n$ implies $f(m)\mid f(n)$ for all $m,n\in N$, then $p\widemid q$ implies $\widetilde{f}(p)\widemid\widetilde{f}(q)$ for all $p,q\in\beta N$.
\end{co}

\dokaz (a) Assume $f(m)\mid m$ for all $m\in N$. Let $A\in \widetilde{f}(p)\cap\cU$. Then $f^{-1}[A]\in p$. By the hypothesis and $A\in\cU$ we have $f^{-1}[A]\ps A$, so $A\in p$ as well. Now $\widetilde{f}(p)\widemid p$ follows from Theorem \ref{Dodp}.

(b) Now assume that $m\mid n$ implies $f(m)\mid f(n)$ for all $m,n\in N$, and that $p\widemid q$. First we note that $A\in\cU$ implies $f^{-1}[A]\in\cU$: if $m\in f^{-1}[A]$ and $m\mid n$, then $f(m)\in A$ and $f(m)\mid f(n)$, so $n\in f^{-1}[A]$ as well. Let $A\in\widetilde{f}(p)\cap\cU$. Then $f^{-1}[A]\in p\cap\cU$, so by $p\widemid q$ we have $f^{-1}[A]\in q\cap\cU$. But this means that $A\in\widetilde{f}(q)$, so by Theorem \ref{Dodp} $\widetilde{f}(p)\widemid\widetilde{f}(q)$.\kraj

We also note that, whenever $f:N\str N$ is a homomorphism, by Proposition \ref{exthom} $p\mid_R q$ implies $\widetilde{f}(p)\mid_R\widetilde{f}(q)$ (and the analogous statements hold for $\mid_L$ and $\mid_M$).

\begin{lm}\label{odnosi}
(a) $D(p)\ps C(p)$ for every $p\in\beta N$.

(b) $p\mid_L q$ implies $p\widemid q$ for all $p,q\in\beta N$.

(c) $p\mid_M q$ implies $p\widemid q$ for all $p,q\in\beta N$.
\end{lm}

\dokaz (a) Suppose $A\in D(p)$. Then the set $B_A=\{m\in N:mN\ps A\}$ is in $p$. But $B_A\ps A/n$ for all $n\in N$, because $mN\ps A$ implies $mn\in A$. Thus $A/n\in p$ for all $n$, so $A\in C(p)$.

(b) From $p\mid_L q$, by Theorem \ref{Codp} we have $C(p)\ps q$. By (a), $D(p)\ps C(p)$, so $D(p)\ps q$. By Theorem \ref{Dodp}, this means that $p\widemid q$.

(c) We note that $p\mid_M q$ means that, for some $r\in\beta N$, $p\mid_L r$ and $r\mid_R q$. Hence $\mid_M$ is the transitive closure of $\mid_L\cup\mid_R$. But, by Theorem \ref{rsttilda} the relation $\widemid$ is transitive, and it contains $\mid_L$ and $\mid_R$, so $\mid_M\ps\widemid$.\kraj

\section{Prime and irreducible elements}

\begin{de}\label{defprime}
An element $p\in\beta N$ is {\it irreducible} in $X\ps\beta N$ if it can not be represented in the form $p=xy$ for $x,y\in X\setminus\{1\}$. $p\in\beta N$ is {\it prime} if $p\mid_R xy$ for $x,y\in\beta N$ implies $p\mid_Rx$ or $p\mid_R y$.
\end{de}

Clearly, $p\in N$ is irreducible (prime) in $N$ if and only if it is irreducible in $\beta N$. We will now see that $p\in N$ is prime in $\beta N$ (by Definition \ref{defprime}) if and only if it is prime in the usual sense.

\begin{lm}
If $n\in N$ is a prime number and $n\mid xy$ for some $x,y\in\beta N$, then $n\mid x$ or $n\mid y$.
\end{lm}

\dokaz $n\mid xy$ means that $nN\in xy$, i.e.\ $\{m\in N:(nN)/m\in y\}\in x$. Now, for $m\in N$ such that $n\nmid m$ we have $(nN)/m=nN$, and for those such that $n\mid m$ we have $(nN)/m=N$. So either $nN\in y$, or $\{m\in N:(nN)/m\in y\}=nN\in x$.\kraj

Irreducible elements of $\beta N$ are exactly the $\mid_L$-minimal and $\mid_R$-minimal (the two notions are the same) ultrafilters. Let $P=\{n\in N:n\mbox{ is prime}\}$.

\begin{lm}\label{irredP}
If $p\in\beta N$ and $P\in p$, then $p$ is irreducible in $\beta N$.
\end{lm}

\dokaz For $p\in N$ the statement is obvious, so let $p\in N^*$. Suppose that $P\in p$, but $p=xy$ for some $x,y\in\beta N$. But $P/n$ is $P$ when $n=1$, $\{1\}$ when $n$ is prime, and empty set if $n$ is composite. So $\{n\in N:P/n\in y\}\in x$ would imply that either $x=1$ or $y=1$.\kraj

\begin{pp}\label{denseirred}
(\cite{HS}, Theorem 6.35) $N^*N^*$ is nowhere dense in $N^*$, i.e.\ for every $A\in[N]^{\aleph_0}$ there is $B\in[A]^{\aleph_0}$ such that all elements of $B^*$ are irreducible in $N^*$.
\end{pp}

Now we can conclude that the reverse of Lemma \ref{irredP} does not hold.

\begin{lm}
There is $p\in \beta N$ irreducible in $\beta N$ such that $P\notin p$.
\end{lm}

\dokaz The set $X=\bigcap_{n\in P}\overline{(nN)^c}\cap\overline{(P\cup\{1\})^c}$ is a nonempty $G_\delta$ set ($X\neq\emptyset$ is shown as in the proof of Theorem \ref{zadatidelioci}). Obviously, it is also a subset of $N^*$. By Proposition \ref{Gdelta} it has nonempty interior: there is $A\in[N]^{\aleph_0}$ such that $A^*\ps X$. By Proposition \ref{denseirred} there is $B\in[A]^{\aleph_0}$ such that all the elements $p\in B^*$ are irreducible in $N^*$. But $p\in\overline{(nN)^c}$ for $n\in P$ implies that also $n\nmid p$ for $n\in N$, so $p$ is actually irreducible in $\beta N$. Finally, since $p\in\overline{P^c}$, it follows that $P\notin p$.\kraj

\begin{lm}\label{nijeRC2}
$\mid_M$ does not satisfy (\ref{RCL}).
\end{lm}

\dokaz By Theorem \ref{zadatidelioci} there is a $G_\delta$ set $X\ps N^*$ such that all its elements are divisible by all elements of $N$. By Proposition \ref{Gdelta} there is a nonempty open set $A^*\ps X$, and by Proposition \ref{denseirred} there is $B^*\ps A^*$ such that all elements of $B^*$ are irreducible. Taking any two different elements $p\in N^*$ and $q\in B^*$ we have that $q$ is divisible by all $n\in N$ but not $\mid_M$-divisible by $p$, an element belonging to the closure of $N$.\kraj

\section{Final words}

The main motivation for the results presented in this paper is to apply them to various number theory problems by translating those problems into $\beta N$. Of course, the method is not convenient for problems of finitary character, even for easy ones. However, many open problems in number theory deal with the existence of infinitely many numbers satisfying certain properties, and it is our hope that translating them into $\beta N$ will give a new insight and help solve some of them.

On the other hand, the section on continuous extensions of binary relations is quite general and perhaps interesting on its own.


\footnotesize
\begin{thebibliography}{99}
\bibitem{HS}
       Hindman N., Strauss D.:
       Algebra in the Stone-\v Cech compactification, theory and applications.
       2nd revised and extended edition, De Gruyter, 2012.
\end{thebibliography}
\end{document}